\documentclass[a4paper,10pt,reqno]{amsart}
\usepackage{amsmath,amssymb}
\usepackage{url}
\usepackage{comment}

\newcommand{\Q}[1]{\mathbb{Q}(\sqrt{#1})}
\newcommand{\Z}{\mathbb{Z}}

\newcommand{\calO}{\mathcal{O}}

\newtheorem{Rmk}{Remark}

\newtheorem{Thm}{Theorem}
\newtheorem{Cor}{Corollary}

\newcommand{\qf}[1]{\langle #1 \rangle}
\newcommand{\conj}[1]{\overline{#1}}
\newcommand{\Case}[1]{\noindent{\bf Case #1.}}
\newcommand{\binlattice}[4]{\begin{pmatrix}
  #1 & #2 \\
  #3 & #4
\end{pmatrix}}
\newcommand{\terlattice}[9]{\begin{pmatrix}
  #1 & #2 & #3 \\
  #4 & #5 & #6 \\
  #7 & #8 & #9
\end{pmatrix}}
\newcommand{\comega}{{\conj\omega}}
\newcommand{\nequiv}{\not\equiv}
\newcommand{\Tr}{\operatorname{Tr}}

\newcommand{\rank}{\operatorname{rank}}
\renewcommand{\Re}{\operatorname{Re}}

\newcommand{\mystrut}{\rule{0pt}{3ex}}

\title[The 15-Theorem for Universal Hermitian Lattices]%
{The Fifteen Theorem for Universal Hermitian Lattices Over
Imaginary Quadratic Fields}%

\keywords{universal Hermitian form}%
\subjclass[2000]{Primary 11E39; Secondary 11E20, 11E41}%

\author[Byeong Moon Kim]{Byeong Moon Kim}
\address{Department of Mathematics, Kangnung National University, Kangnung, Korea}%
\email{kbm@kangnung.ac.kr}

\author[Ji Young Kim]{Ji Young Kim}
\address{School of Mathematics, Korea Institute for Advanced Study, Hoegiro 87, Dongdaemun-gu, Seoul, 130-722, Korea.}%
\email{jykim@kias.re.kr}

\author[Poo-Sung Park]{Poo-Sung Park}
\address{Korea Institute for Advanced Study, Cheongnyangni 2-dong, Dongdaemun-gu, Seoul, 130-722, Korea}%
\email{sung@kias.re.kr}

\thanks{The second and third authors were partially supported by KRF(2003-070-c00001)}

\begin{document}

\begin{abstract}
We will introduce a method to get all universal Hermitian lattices over imaginary quadratic fields
$\Q{-m}$ for all $m$. For each imaginary quadratic field $\Q{-m}$, we obtain a criterion on
universality of Hermitian lattices: if a Hermitian lattice $L$ represents 1, 2, 3, 5, 6, 7, 10, 13,
14 and 15, then $L$ is universal. We call this the \emph{fifteen theorem for universal Hermitian
lattices}. Note that the difference between Conway-Schneeberger's fifteen theorem and ours is the
number 13.
\end{abstract}

\maketitle

\section{introduction}

The research for positive definite rational quadratic forms for which the represented integer set
is as large as possible has its origins at the beginning of modern number theory. In 1770, Lagrange
\cite{jlL_70} found the famous \emph{four square theorem}: the positive definite quadratic form
$x_1^2 + x_2^2 + x_3^2 +x_4^2$ represents all positive integers. Since then, his theorem has been
generalized in many directions. One of the generalizations is to find all positive definite
quadratic forms that represent all positive integers, which we call universal quadratic forms. The
first breakthrough in this direction was made by Ramanujan \cite{sR_17}. In 1917, he discovered all
55 quaternary diagonal universal forms, up to isometry. In 1927, Dickson \cite{leD_27} confirmed
Ramanujan's list except one form which was not universal and extended Ramanujan's results to
non-diagonal forms. It was Dickson who called those forms universal. In 1948, Willerding
\cite{mfW_48} found 124 quaternary \emph{classical} non-diagonal universal forms, up to isometry,
and claimed that the list was complete. But her list was incomplete with some mistakes.

More generally, a positive definite quadratic form over a totally real number field is said to be
universal if every totally positive integer of the field is represented by it. To establish a
general context for the following discussion, let $F$ be a totally real algebraic number field with
ring of integers $\calO$ and let $\calO^{+}$ be the subring of totally positive elements of
$\calO$. By the term \emph{ $\calO$-lattice} (or, \emph{lattice over $F$}), we will always mean an
$\calO$-module on the totally positive definite quadratic space $(V,Q)$ over $F$. As the similar
way of defining an quadratic forms (over $\Z$), we can define an \emph{quadratic lattices over
$\calO$}. So, an $\calO$-lattice $L$ is said to be universal if every element of $\calO^{+}$ can be
represented by $L$. In 1928, G\"otzky \cite{fG_28} proved that $x_1^2 + x_2^2 + x_3^2 + x_4^2$ is
universal over $\mathbb{Q}(\sqrt{5})$. In 1941, Maass \cite{hM_41} proved the \emph{three square
theorem}, which states: the quadratic form $x_1^2 + x_2^2 + x_3^2$ is universal over $\Q{5}$. All
positive definite ternary universal forms over real quadratic fields were determined in
\cite{wkC_mhK_sR_96}. Further developments on universal forms over totally real number fields are
established by B. M. Kim (see \cite{bmK1}, \cite{bmK_99} and \cite{bmK_00}).

If a Hermitian lattice represents all positive integers, we call it universal. In 1997, Earnest and
Khosravani \cite{agE_aK_97} found 13 universal binary Hermitian forms over imaginary quadratic
fields of class number 1. If the quadratic field over $\mathbb{Q}$ has the class number bigger than
1, Iwabuchi \cite{hI_00} determined all universal binary Hermitian lattices over this field. After
that, J.-H. Kim and P.-S. Park \cite{jhK_psP} added three binary Hermitian forms to the
Earnest-Khosravani-Iwabuchi's list and completed the list. Further generalizations of Hermitian
lattices were made by P.-S. Park \cite{psP_05} and A. Rokicki \cite{aR_05}.

In 1997, Conway and Schneeberger announced the \emph{fifteen
theorem} for classical universal quadratic forms, which
characterizes the universality by representability of a finite set
of numbers, namely, 1, 2, 3, 5, 6, 7, 10, 14, and 15 (see
\cite{jhC_wS}). Using the fifteen theorem, they corrected several
mistakes in Willerding's list and announced the new and complete
list of 204 quaternary classical universal forms up to isometry. B.
M. Kim, M.-H. Kim and B.-K. Oh proved 2-universality analogy in
\cite{bmK_mhK_bkO_99}. Recently, Bhargava and Hanke enunciated that
they proved the \emph{290-conjecture} which characterizes the
universality of (nonclassical) quadratic forms. That is, if a
(nonclassical) quadratic form represents 29 numbers, say, 1, 2, 3,
5, 6, 7, 10, 13, 14, 15, 17, 19, 21, 22, 23, 26, 29, 30, 31, 34, 35,
37, 42, 58, 93, 110, 145, 203, 290, then it is universal (see
\cite{jH_06}). One can use this theorem to prove the universality of
a Hermitian lattice. But we succeeded in proving the universality of
a Hermitian lattice without this big theorem.

In \cite{jhC_00}, Bhargava's generalization of the fifteen theorem was announced: for any infinite
set $S$ of positive integers, there is a finite subset $S_0$ of $S$ such that any positive definite
quadratic form that represents every element of $S_0$ represents all elements of $S$. And he found
$S_0$ for some interesting sets $S$. In \cite{bmK_mhK_bkO_05}, B. M. Kim, M.-H. Kim and B.-K. Oh
proved the \emph{finiteness theorem} as a generalization of Bhargava's result: for any infinite set
$S$ of positive definite quadratic forms of bounded rank, there is a finite subset $S_0$ of $S$
such that any positive definite quadratic form that represents every element of $S_0$ represents
all of $S$. Nice survey papers related to these subjects are \cite{mhK} and \cite{jH_04}.

In this paper, first, we will suggest a matrix representation for
non-free Hermitian lattices. Due to this matrix representation, we
can do escalate to find candidates of universal Hermitian lattices
including non-free lattices over imaginary quadratic fields
$\Q{-m}$. Second, we obtain a Conway-Schneeberger-Bhargava type
criterion on universality of Hermitian lattices: if a Hermitian
lattice $L$ represents 1, 2, 3, 5, 6, 7, 10, 13, 14 and 15, then $L$
is universal. Hence we call this theorem the \emph{fifteen theorem
for universal Hermitian lattices}. As the language of finiteness
theorem for representability, $S$ is the set of all positive
integers and $S_0 = \{1, 2, 3, 5, 6, 7, 10, 13, 14, 15 \}$. We will
call the set $S_0$ a \emph{set of critical numbers}. Even though the
fifteen theorem and 290-theorem gives us the rough upper bound for
the critical numbers that lie in between 15 to 290, it is hard to
figure out the set of critical numbers for each field. For each
imaginary quadratic field $\Q{-m}$,  we will give an optimal set of
critical numbers by arithmetic calculation. For example, a Hermitian
lattice $L$ over $\Q{-39}$ is universal if and only if $L$
represents 1, 2, 3, 5, 6, 7, 13.

\section{Preliminaries}

The notation and terminology of O'Meara's book \cite{oO_73} will be adopted here. For the
terminology specific to the Hermitian case, the paper \cite{nJ_40} can be referred to.  We begin by
setting some additional notations that will remain in effect throughout this paper. Let $F$ denote
the imaginary quadratic field $\Q{-m}$ for a positive square-free integer $m$ with nontrivial
$\mathbb{Q}$-involution and let $\calO$ be the ring of integers of $F$. It is well-known that
$\calO$ is generated by $\{1, \omega=\omega_{m}\}$ over $\Z$, where $\omega_{m} = \sqrt{-m}$ if $m
\equiv 1, 2 \pmod{4}$ or $\omega_{m} = \frac{1+\sqrt{-m}}2$ if $m \equiv 3 \pmod{4}$. By the term
\emph{ $\calO$-lattice} $L$ (or \emph{integral lattice $L$ over $F$}), we will mean a finitely
generated $\calO$-module on the Hermitian space $(V,H)$ over $F$, where $V$ is an $n$-dimensional
vector space over $F$ with the nondegenerate Hermitian form $H$. All lattices considered here will
be assumed to be integral and positive definite in the sense that $H(x,y) \in \calO$ for all $x, y
\in L$ and $H(x):=H(x,x) > 0$ for all $x \ne 0$. It follows from these assumptions that $H(x)$,
called a (Hermitian) norm, is in $\mathbb{Z}$ for all $x \in L$.

As the ring of integers of an imaginary quadratic field is not generally a principal ideal domain,
lattices do not need to be free. Let $\{ v_1, \dotsc, v_n \}$ be an $\calO$-basis for $L$. In case
that $\calO$ is a principal ideal domain, every Hermitian lattice is free. Therefore $L = \calO v_1
+ \dotsb + \calO v_n$ and there is a function $f : \calO^n \longrightarrow \mathbb{Z}$ defined by
$f(x_1, \dotsc, x_n) = H(\sum x_i v_i) = \sum H(v_i, v_j) x_i \conj{x_j}$. Such a function will be
referred to as a Hermitian form associated to $L$. And we can obtain an associated Hermitian matrix
for $L$ by taking the $n \times n$-matrix whose entry is $H(v_i, v_j)$. If a basis $\{v_1 \dotsc,
v_n \}$ for $L$ is orthogonal then the associated matrix of $L$ is denoted by $\qf{H(v_1), \dotsc,
H(v_n)}$. Similarly, in case that $\calO$ is not a principal ideal domain, there is a fractional
ideal $\mathfrak{A}$ such that $L= \calO v_1 + \dotsb + \calO v_{n-1} + \mathfrak{A} v_n$ by
\cite[81:5]{oO_73}. Since any ideal in $\calO$ is generated by at most two elements, we can write
$L= \calO v_1 + \dotsb + \calO v_{n-1} + (\alpha, \beta)\calO v_n$ for some $\alpha, \beta \in
\calO$. Therefore, we have a Hermitian form $f : \calO^{n-1} \times \calO \longrightarrow
\mathbb{Z}$ associated to $L$ such that $f(x_1, \dotsc, x_n) = H(\sum _{i=1} ^{n-1} x_i v_i + x_n
\alpha v_n + x_n \beta v_n)$. Also we have a Hermitian $(n+1) \times (n+1)$-matrix associated to
$L$ as follows:
\[
    \begin{pmatrix}
        H(v_1, v_1)         & \dotsc & H(v_1, \alpha v_n)        & H(v_1, \beta v_n)\\
        \vdots              & \ddots & \vdots                    &\vdots\\
        H(\alpha v_n,  v_1) & \dotsc & H(\alpha v_n, \alpha v_n) & H(\alpha v_n, \beta v_n)\\
        H(\beta  v_n,  v_1) & \dotsc & H(\beta  v_n, \alpha  v_n) & H(\beta  v_n, \beta v_n)\\
    \end{pmatrix}.
\]
Note that this matrix is positive semi-definite, but this represents an $n$-ary positive definite
Hermitian lattice.

Considering a $2n$-dimensional vector space $\widetilde{V}$ over $\mathbb{Q}$ corresponding to $V$
as defined in \cite{nJ_40}, we can regard $(V,H)$ over $F$ as a $2n$-dimensional quadratic space
$(\widetilde{V}, B_H)$ such that $B_H (x,y) = \frac{1}{2} [ H(x,y) + H(y,x) ] = \frac{1}{2}
\Tr_{F/\mathbb{Q}}(H(x,y))$. Analogously, by viewing $L$ as a $\mathbb{Z}$-lattice on
$(\widetilde{V}, B_H)$ we can obtain a quadratic $\mathbb{Z}$-lattice $\widetilde{L}$ on
$\widetilde{V}$ associated to a Hermitian $\calO$-lattice $L$ on $V$ and also $\widetilde{f}(x_1,
y_1, \dotsc, x_n, y_n) = f(x_1 +\omega y_1, \dotsc, x_n +\omega y_n)$ is an associated quadratic
form in $2n$-variables corresponding to this lattice $\widetilde{L}$. For the convenience, we say
that $\widetilde{f}$ is an associated quadratic form of $L$. For example, the associated quadratic
form of the Hermitian lattice $\qf{1}$ over $\Q{-m}$ is
\[
    \begin{cases}
        ~x_1^2 + m y_1^2                       &\text{ if } m \equiv 1, 2 \pmod{4}\\
        ~x_1^2 + x_1 y_1 + \frac{1+m}{4} y_1^2 &\text{ if } m \equiv 3 \hspace{2.2ex}   \pmod{4}
    \end{cases}
\]
To distinguish the matrix of a quadratic $\Z$-lattice $\widetilde{L}$ from the matrix of a
Hermitian lattice $L$, we will add subscript $\Z$ to the matrix of the quadratic $\Z$-lattice
$\widetilde{L}$. And a Hermitian $\calO$-lattice $L$ represents a $\Z$-lattice $M$ if and only if
the associated quadratic form of $L$ represents $M$.

If $L$ is a universal lattice of rank $n$, then there are infinitely many universal lattices of
rank $k$ ($k > n$) which contain $L$. Thus in order to obtain any meaningful finiteness result for
such lattices, we should consider a \emph{new} universal lattice. A universal lattice is called new
when it does not contain any universal lattice of smaller rank.

For $m \equiv 1, 2 \pmod{4}$, the associated quadratic lattice of a Hermitian lattice over $\Q{-m}$
is a classical $\Z$-lattice. Thus we can determine the universality of a Hermitian lattice over
$\Q{-m}$ via applying the fifteen theorem. On the contrary, for $m \equiv 3 \pmod{4}$, the
associated quadratic lattice of a Hermitian lattice over $\Q{-m}$ is a nonclassical $\Z$-lattice.
Therefore, the recent big result of Bhargava and Hanke, the 290-conjecture, can be applied to prove
the universality.

We adopt some notations from Conway-Sloane \cite{jhC_00}. The notation $p^{d}$ (resp. $p^{e}$)
denotes an odd (resp. even) power of $p$; if $p=2$, $u_k$ denotes a unit of the form $8n + k$
($k=1,3,5,7$) and if $p$ is odd, $u_{+}$ (resp. $u_{-}$) denotes a unit which is a quadratic
residue (resp. non-residue) modulo $p$.

\section{main results}
Note that $x\conj{x}+y\conj{y}+z\conj{z}+u\conj{u}$ represents all
positive integers even when all variables take values in $\Z$. Thus
it is obviously a quaternary universal Hermitian lattice over all
imaginary quadratic fields. Since classical universal quadratic
forms over $\Z$ was already classified, we do not need to
investigate this kind of lattices. A Hermitian lattice is called
\emph{inherited} if its coefficients are all rational integers. If a
Hermitian lattice is not inherited, we call it \emph{uninherited}.
We will mainly consider uninherited universal Hermitian lattices.

If a lattice $L$ is not universal, define the \emph{truant} of $L$
to be the smallest positive integer not represented by $L$. An
escalation of a nontrivial lattice $L$ is defined to be any lattice
which is generated by $L$ and a vector whose norm is equal to the
truant of $L$. Conway, Schneeberger and Bhargava suggested this
escalation method for free lattices (see \cite{jhC_wS},
\cite{mB_00}). In this article, we suggested the method for a matrix
representation for non-free Hermitian lattices (see section 2) and
we use the escalation method to find candidates of universal
Hermitian lattices including non-free Hermitian lattices.

Let $L$ be a universal Hermitian lattice over $\Q{-m}$. Since $L$ represents $1$, $L \cong \qf{1}
\perp M$ for some lattice $M$. Since $2 \to L$, there are three incompatible cases: $2 \to \qf{1}$,
$1 \to M$, or $2 \to M$.

\Case{I} $2 \to \qf{1}$: In this case, $\alpha \conj{\alpha} = 2$
has a solution for some $\alpha \in \calO$. Let $\alpha = x_1 + x_2
\omega_m$. If $m \equiv 1, 2 \pmod{4}$, then $\alpha \conj{\alpha} =
x_1^2 + m x_2^2 =2$. So $m = 1$ or $2$. If $m \equiv 3 \pmod{4}$,
then $\alpha \conj{\alpha} = x_1^2 + x_1 x_2 + \frac{m+1}{4} x_2^2
=2$. So $m = 7$. Therefore all possible $m$'s are $1,2;7$. We use a
semicolon to distinguish $m \equiv 1, 2 \pmod{4}$ and $m \equiv 3
\pmod{4}$. To represent from 1 through 5, $L$ should have the
following sublattices:

\begin{tabular}{ll}
$\Q{-1}$: & $\qf{1,1}, \qf{1,2}, \qf{1,3}$, \\
$\Q{-2}$: & $\qf{1,1}, \qf{1,2}, \qf{1,3}, \qf{1,4}, \qf{1,5}$, \\
$\Q{-7}$: & $\qf{1,1}, \qf{1,2}, \qf{1,3}$.
\end{tabular}

These Hermitian lattices are all universal by \cite{agE_aK_97} and \cite{jhK_psP} and we know that
there are no other new Hermitian universal lattices over $\Q{-m}$ for $m=1,2;7$. Hence we assume
that $m \ne 1,2;7$.

\Case{II} $2 \not\to \qf{1}$ and $1 \to M$: Then $L \cong \qf{1,1} \perp K$ for some lattice $K$.
Since $3$ is a truant of $\Z$-lattice $\qf{1,1}_\Z$ and $L$ should represent $3$, we have four
cases.

\Case{II-1} $3 \to \qf{1,1}$: In this case, $\alpha \conj{\alpha} +
\beta \conj{\beta} = 3$ has a solution for some $\alpha$, $\beta \in
\calO$, the possible $m$'s are $3$ or $11$. The universalities of
$\qf{1,1}$ over $\Q{-3}$ and $\Q{-11}$ were shown in
\cite{agE_aK_97}. Assume that $m \ne 3,11$ for the following
subcases.

\Case{II-2} $1 \to K$: Then $L \cong \qf{1,1,1} \perp N$ for some lattice $N$. Note that the truant
of $\qf{1,1,1}_\Z$ is 7.

If $7 \to \qf{1,1,1}$ over $\Q{-m}$, then we have $m = 5,6;15,19,23$. Since the sublattice
$\qf{1,2}$ of $\qf{1,1,1}$ is universal over $\Q{-m}$ with $m=5;19$ by \cite{agE_aK_97} and
\cite{hI_00}, we do not need to consider the cases $m=5;19$.

\Case{II-2(1)} When $m=6$, the lattice $\qf{1,1,1}$ represents a universal quadratic form
$\qf{1,1,1,6}_\Z$, i.e., the associated quadratic form $\qf{1,1,1,6,6,6}_\Z$ of $\qf{1,1,1}$
represents a universal quadratic form $\qf{1,1,1,6}_\Z$. Thus the Hermitian lattice $\qf{1,1,1}$ is
also universal.

\Case{II-2(2)} When $m=15,23$, $\qf{1,1,1}$ can be written as an associated quadratic form over
$\Z$:
\[
    x_1^2+x_1x_2+\frac{1+m}4x_2^2 + y_1^2+y_1y_2+\frac{1+m}4y_2^2 +
    z_1^2+z_1z_2+\frac{1+m}4z_2^2.
\]
This quadratic form has a sublattice $\qf{1,1,1}_\Z$. It represents all positive integers except
the form $2^{e} u_7$. If we set $x_2=2$ and $y_2 = z_2=0$, then the quadratic form becomes
\[
    (x_1+1)^2 + y_1^2 + z_1^2 + m.
\]
And if we set $x_2=y_2=2$ and $z_2=0$, then the quadratic form becomes
\[
    (x_1+1)^2 + (y_1+1)^2 + z_1^2 + 2m.
\]

Let $n=2^{e} u_7$. If $n \ge 2m$ then at least one of $n-m$ and $n-2m$ is not of the form $2^{e}
u_7$. Thus one of $n-m$ or $n-2m$ can be represented by $\qf{1,1,1}_\Z$. That is, $n \to
\qf{1,1,1}$ over $\Q{-m}$. It is easy to verify that $n \rightarrow \qf{1,1,1}$ for $n < 2m$.

Hence we obtain new universal Hermitian lattices $\qf{1,1,1}$ over $\Q{-m}$ with $m = 6;15,23$.

If $m \ne 5,6;15,19,23$, $\qf{1,1,1}$ cannot represent 7 over $\Q{-m}$ and thus $N$ should
represent 1, 2, 3, 4, 5, 6 or 7. But the quadratic form $\qf{1,1,1,a}_\Z$ is universal for $a = 1,
2 \dotsc, 7$ by \cite{sR_17}. That is, these lattices $\qf{1,1,1,a}$ are \emph{inheritedly}
universal when $a=1, 2, \dotsc, 7$.

\Case{II-3} $1 \not\to K$ and $2 \to K$: Then $L$ contains a ternary lattice $\qf{1,1,2}$. Note
that the truant of $\qf{1,1,2}_\Z$ is 14.

\Case{II-3 a} $14 \to \qf{1,1,2}$: Then we have
\[
    m =  5, 6, 10, 13, 14; 15, 19, 23, 31, 35, 39, 43, 47, 51, 55.
\]

Since the sublattice $\qf{1,2}$ of $\qf{1,1,2}$ is universal over $\Q{-5}$ and $\Q{-19}$,
$\qf{1,1,2}$ is not a new universal Hermitian form over $\Q{-5}$ and $\Q{-19}$.

\Case{II-3 a(1)} $m \nequiv 3 \pmod{4}$: Since $m=6, 10, 13, 14$, the Hermitian lattice
$\qf{1,1,2}$ over $\Q{-m}$ represents a universal quadratic lattice $\qf{1,1,2,m}_\Z$. Therefore
$\qf{1,1,2}$ is universal over $\Q{-m}$.

\Case{II-3 a(2)} $m \equiv 3 \pmod{4}$: $m = 15$, $23$, $31$, $35$, $39$, $43$, $47$, $51$, $55$.
For these cases we can write $\qf{1,1,2}$ as an associated quadratic form over $\Z$:
\[
    x_1^2+x_1x_2+\frac{1+m}4x_2^2 + y_1^2+y_1y_2+\frac{1+m}4y_2^2 +
    2\left(z_1^2+z_1z_2+\frac{1+m}4z_2^2\right).
\]
This form has a sublattice $\qf{1,1,2}_\Z$ and it represents all positive integers except the form
$2^{d}u_7$. If we set $x_2=2$ and $y_2 = z_2=0$, then the quadratic form becomes
\[
    (x_1+1)^2 + y_1^2 + 2z_1^2 + m.
\]

Let  $n=2^{d} u_7$. If $n \ge m$, then $n-m$ is not of the form $2^{d}u_7$ since $m$ is odd. Thus
$n \to \qf{1,1,2}$ over $\Q{-m}$. And we can show that $n \to \qf{1,1,2}$ for $n < m$ by direct
calculation.

Hence $\qf{1,1,2}$ is a new universal Hermitian lattice over $\Q{-m}$ for
\[
    m = 6, 10, 13, 14; 15, 23, 31, 35, 39, 43, 47, 51, 55.
\]

\Case{II-3 b} $14 \not\to \qf{1,1,2}$: From $14 \not\to \qf{1,1,2}$, we have
\[
    \begin{cases}
        ~m \ge 17 ~\text{ if } ~m \nequiv 3 \pmod{4},\\
        ~m \ge 59 ~\text{ if } ~m \equiv  3 \pmod{4}.
    \end{cases}
\]
Since $14 \to L$, $L$ should contain an escalation lattice $\ell$ of $\qf{1,1,2}$ as follows:
\[  \ell \cong
    \begin{pmatrix}
        1 & 0 & 0 & \alpha \\
        0 & 1 & 0 & \beta \\
        0 & 0 & 2 & \gamma \\
        \conj{\alpha} & \conj{\beta} & \conj{\gamma} & 14
    \end{pmatrix}
    \cong
    \begin{pmatrix}
        1 & 0 & 0 & 0 \\
        0 & 1 & 0 & 0 \\
        0 & 0 & 2 & \gamma \\
        0 & 0 & \conj{\gamma} & 14-\alpha\conj{\alpha}-\beta\conj{\beta}
    \end{pmatrix}
\]
for some $\alpha$, $\beta$, $\gamma$. Note that $14 -\alpha \conj{\alpha} - \beta \conj{\beta} = 4,
5, 6, 9, 10, 12, 13, 14$. The above lattice $\ell$ can be reduced once more to $\qf{1,1} \perp
\binlattice2a{\conj{a}}b$ with $a = 0, 1, \omega, -1+\omega$ and suitable $b$. If $a=0$ or $1$,
then $\ell$ is one of the followings:
\begin{align*}
    &\qf{1,1}\perp\binlattice200b \text{ with } b=2,\dotsc,14, \\
    &\qf{1,1}\perp\binlattice211b \text{ with } b=2,4,5,6,8,9,10,12,13,14.
\end{align*}
These are all inherited universal lattices for all the above $m$'s. Now consider the case of
$a=\omega$ or $-1+\omega$.

\Case{II-3 b(1)} $m \nequiv 3 \pmod{4}$: From the positive
semi-definiteness of $\ell$, we have $m=17,21,22,26$. In addition,
\begin{align*}
    &\text{ if } a=\omega, \text{ then }  b = \begin{cases}
                                             9,  10, 12, 13, 14,  &\text{ when } m = 17,\\
                                             12, 13, 14,          &\text{ when } m = 21, 22,\\
                                             13, 14,                  &\text{ when } m = 26,\\
                                            \end{cases}\\
    &\text{ and if } a= -1+\omega, \text{ then }  b = \begin{cases}
                                                 9,  10, 12, 13, 14, &\text{ when } m = 17,\\
                                                 12, 13, 14,         &\text{ when } m = 21, 22,\\
                                                 14,                     &\text{ when } m = 26.
                                                 \end{cases}
\end{align*}
Since each $\ell$ represents a universal quaternary quadratic lattice $\qf{1,1,2,b}_\Z$ or
$\qf{1,1}_\Z \perp\binlattice211b_\Z$ for each $b$, it is universal. Note that the following
lattices are the first examples of non-free ternary lattices:
\[
    \begin{array}{ll}
        \qf{1,1} \perp \binlattice2{-1+\omega_{17}}{-1+\comega_{17}}9,
       &\qf{1,1} \perp \binlattice2{\omega_{26}}{\comega_{26}}{13}.
    \end{array}
\]

\Case{II-3 b(2)} $m \equiv 3 \pmod{4}$: From the positive semi-definiteness of $\ell$, we have $m =
59$, $67$, $71$, $79$, $83$, $87$, $91$, $95$, $103$, $107$, $111$. If $\ell$ represents all
positive integers smaller than $m$, then $\ell$ is universal, because $\ell$ contains $\qf{1,1,2}$.
We can easily check it by direct calculation.

\Case{II-4} $1, 2 \not\to K$ and $3 \to K$: Then $L$ contains a ternary lattice $\qf{1,1,3}$. Note
that the truant of $\qf{1,1,3}_\Z$ is 6.

\Case{II-4 a} $6 \to \qf{1,1,3}$: Then we have
\[
    m = 5, 6; 15, 19, 23.
\]

\Case{II-4 a(1)} $m \nequiv 3 \pmod{4}$: Since $m = 5$ or $6$, the Hermitian lattice $\qf{1,1,3}$
represents a universal quadratic form $\qf{1,1,3,5}_\Z$ or $\qf{1,1,3,6}_\Z$. Thus $\qf{1,1,3}$ is
universal over $\Q{-m}$.

\Case{II-4 a(2)} $m \equiv 3 \pmod{4}$: We have that $\qf{1,1,3}_\Z$ can represent all positive
integers except the form $3^{d} u_{-}$. Let $n=3^{d}u_{-}$.

First consider the case of $m=15$. The associated quadratic form
of $\qf{1,1,3}$ is a nonclassical quadratic form
\[
    x_1^2 + x_1x_2 + 4x_2^2 + y_1^2 + y_1y_2 + 4y_2^2 + 3z_1^2+3z_1z_2+12z_2^2.
\]
Note that if $x_2 = 2$ and $y_2 = z_2 = 0$, then the form becomes
\[
    (x_1+1)^2 + y_1^2 + 3z_1^2 + 15,
\]
and if $x_2 = y_2 = 2$ and $z_2 = 0$, then the form becomes
\[
    (x_1+1)^2 + (y_1+1)^2 + 3z_1^2 + 30.
\]
If $n \ge 30$ then at least one of  $n-15$ and $n-30$ is not of the form $3^{d}u_{-}$. That is, $n$
can be represented by $\qf{1,1,3}$. It is easily verified that $\qf{1,1,3}$ represents all positive
integers less than $30$. Thus it is a new universal Hermitian lattice over $\Q{-15}$.

Similarly, it can be shown that $\qf{1,1,3}$ is a new universal Hermitian lattice over $\Q{-m}$
with $m = 19, 23$ by checking $n-19$ and $n-23$.

\Case{II-4 b} $6 \not\to \qf{1,1,3}$: From $6 \not\to \qf{1,1,3}$,
we have
\[
    \begin{cases}
        ~m \ge 10 ~\text{ if } ~m \nequiv 3 \pmod{4},\\
        ~m \ge 31 ~\text{ if } ~m \equiv  3 \pmod{4}.
    \end{cases}
\]
Since $6 \to L$, $L$ should contain an escalation lattice $\ell$
of $\qf{1,1,3}$ as follows:
\[  \ell \cong
    \begin{pmatrix}
        1 & 0 & 0 & \alpha \\
        0 & 1 & 0 & \beta \\
        0 & 0 & 3 & \gamma \\
        \conj{\alpha} & \conj{\beta} & \conj{\gamma} & 6
    \end{pmatrix}
    \cong
    \begin{pmatrix}
        1 & 0 & 0 & 0 \\
        0 & 1 & 0 & 0 \\
        0 & 0 & 3 & \gamma \\
        0 & 0 & \conj{\gamma} & 6-\alpha\conj{\alpha}-\beta\conj{\beta}
    \end{pmatrix}
\]
for some $\alpha$, $\beta$, $\gamma$. Note that $6 - \alpha
\conj{\alpha} - \beta \conj{\beta} = 4, 5, 6$. Hence the above
lattice $\ell$ can be reduced to $\qf{1,1} \perp
\binlattice3a{\conj{a}}b$ with $a = 0, 1, \omega, 1+\omega,
-1+\omega$ and suitable $b$. If $a=0$ or $1$, then $b = 3, 4, 5, 6$
and $\ell$ is inherited and universal over $\Q{-m}$ for all the
above $m$. Now assume that $a= \omega$, $1+\omega$ or $-1+\omega$.

\Case{II-4 b(1)} $m \nequiv 3 \pmod{4}$: From the positive semi-definiteness of $\ell$, we have
\[
    m = 10, 13, 14, 17.
\]
In addition, if $a= \omega$, $1+\omega$ or $-1+\omega$, then
\begin{align*}
     b = \begin{cases}
         4, 5, 6,  &\text{ when } m = 10,\\
         5, 6,     &\text{ when } m = 13, 14,\\
         6,        &\text{ when } m = 17.\\
         \end{cases}
\end{align*}
Since each $\ell$ represents a universal quaternary quadratic lattice $\qf{1,1,3,b}_\Z$ or
$\qf{1,1}_\Z\perp\binlattice311b_\Z$ for each $b$, it is a universal Hermitian lattice.

\Case{II-4 b(2)} $m \equiv 3 \pmod{4}$: From the positive semi-definiteness of $\ell$, we have $m =
31, 35, 39, 43, 47, 51, 55, 59, 67, 71$. We will apply the same argument in Case II-4 a(2) for
these $m$'s. If $m = 39,51$, then it is enough to check that each $\ell$ represents all positive
integers smaller than $2m$. If $m = 31, 35, 43, 47, 55, 59, 67, 71$, then it is enough to check
whether $\ell$  represent all positive integers smaller than $m$. We can easily check it by direct
calculation, hence $\ell$ is a universal Hermitian lattice.

\Case{III} $2 \not\to \qf{1}$, $1 \not\to M$ and $2 \to M$: Then $L$ contains a binary lattice
$\qf{1,2}$. Note that the truant of $\qf{1,2}_\Z$ is 5. When $m = 5; 3, 11, 19$, $\qf{1,2}$
represents $5$ and $\qf{1,2}$ is universal for these $m$'s by \cite{agE_aK_97} and \cite{hI_00}.
Thus through the Case III, we may assume $m \ne 5; 3, 11, 19$. Since $5 \to L$, the escalation
lattice of $\qf{1,2}$ is
\[
    \terlattice10\alpha02\beta{\conj{\alpha}}{\conj{\beta}}5
    \cong
    \terlattice10002\beta0{\conj{\beta}}{5-\alpha\conj{\alpha}}
\]
for some $\alpha$, $\beta$. This lattice can be reduced to $\qf{1} \perp\binlattice2a{\conj{a}}b$
with $a = 0, 1, \omega$ or $-1+\omega$ and suitable $b$. If $a = 0, 1$, then the escalation
lattices are
\[
    \qf{1,2,2},
    \qf{1,2,3},
    \qf{1,2,4},
    \qf{1,2,5},
    \qf{1} \perp \binlattice2114 \text{ and }
    \qf{1} \perp \binlattice2115
\]
for all the above $m$'s. And if $a=\omega$ or $-1+\omega$, then the escalation lattice $\qf{1}
\perp\binlattice2a{\conj{a}}b$ over $\Q{-m}$ satisfies the following conditions, up to isometry:

\begin{table}[!h]
\begin{tabular}{l|l|l}\hline
$\Q{-m}$  & $a$              & $b$ \\ \hline
$\Q{-6}$  & $\omega_{6}$     & $3, 4, 5$, \\
$\null$   & $-1+\omega_{6}$  & $4, 5$, \\
$\Q{-10}$ & $\omega_{10}$    & $5$, \\
$\Q{-15}$ & $\omega_{15}$    & $2, 3, 4, 5$, \\
$\null$   & $-1+\omega_{15}$ & $3, 4, 5$, \\
$\Q{-23}$ & $\omega_{23}$    & $3, 4, 5$, \\
$\null$   & $-1+\omega_{23}$ & $3, 4, 5$, \\
$\Q{-31}$ & $\omega_{31}$    & $4, 5$,\\
$\null$   & $-1+\omega_{31}$ & $4, 5$, \\
$\Q{-35}$ & $\omega_{35}$    & $5$, \\
$\null$   & $-1+\omega_{35}$ & $5,$ \\
$\Q{-39}$ & $\omega_{39}$    & $5$, \\
$\null$   & $-1+\omega_{39}$ & $5$. \\ \hline
\end{tabular} 
\caption{Conditions for {\small$\protect\qf{1}\perp\protect\binlattice2a{\protect\conj{a}}b$}}%
\label{tbl:escalation_lattices_in_Case_III}%
\end{table}

Now, we will treat each lattice in Case III-1 to Case III-7.

\Case{III-1} $\qf{1,2,2} \to L$: Note that $\qf{1,2,2}$ represents a quadratic lattice
$\qf{1,2,2}_\Z$ whose truant is $7$.

\Case{III-1 a} $7 \to \qf{1,2,2}$: Then we have $m=6$ and $\qf{1,2,2}$ is a new universal Hermitian
lattice as $\qf{1,2,2}$ represents a universal quadratic lattice $\qf{1,2,2,6}_\Z$. Through the
Case III-1, we may assume $m \ne 6$.

\Case{III-1 b} $7 \not\to \qf{1,2,2}$: Then $L$ should contain an escalation lattice $\ell$ of
$\qf{1,2,2}$ as following,
\[
    \ell \cong
    \begin{pmatrix}
        1 & 0 & 0 & \alpha \\
        0 & 2 & 0 & \beta \\
        0 & 0 & 2 & \gamma \\
        \conj{\alpha} & \conj{\beta} & \conj{\gamma} & 7
    \end{pmatrix}
\]
for some $\alpha, \beta, \gamma$.

\Case{III-1 b(1)} $m \nequiv 3 \pmod 4$: Then $\ell$ represents a quadratic lattice
\[
    \ell' = \begin{pmatrix}
        1 & 0 & 0 & \Re{\alpha} \\
        0 & 2 & 0 & \Re{\beta} \\
        0 & 0 & 2 & \Re{\gamma} \\
        \Re{\alpha} & \Re{\beta} & \Re{\gamma} & 7
    \end{pmatrix}_\Z.
\]
There are 16 lattices of above form up to isometry and they are all universal by the fifteen
theorem.

\Case{III-1 b(2)} $m \equiv 3 \pmod{4}$: From the positive semi-definiteness, $\ell$'s are all
inherited for $m \ge 59$. Then, since $\ell$ represents the above quadratic lattice  $\ell'$,
$\ell$ is universal. We may assume that $m=15$, $23$, $31$, $35$, $39$, $43$, $47$, $51$, $55$.
Note that $\qf{1,2,2}_\Z$ represents all positive integers only except the form $2^{e} u_7$ and
$\qf{1,2,2}$ represents $\qf{1,2,2,m}_\Z$.

Let $n=2^{e} u_7$. If $n \ge 4m$, then at least one of  $n-m$ and $n-4m$ is not of the form $2^{e}
u_7$. Thus $n \to \qf{1,2,2}_\Z$. If $n < 4m$, then we can show that $n$ is represented by each
$\ell$ by direct calculation. Hence the universality of $\ell$ is proved.

\Case{III-2} $\qf{1,2,3} \to L$: Note that $\qf{1,2,3}$ represents a quadratic lattice
$\qf{1,2,3}_\Z$ whose truant is $10$.

\Case{III-2 a} $10 \to \qf{1,2,3}$: Then we have $m=6$, $10$; $15$, $23$, $31$, $39$.

\Case{III-2 a(1)} $m \nequiv 3 \pmod{4}$: Since $m = 6$ or $10$, the Hermitian lattice $\qf{1,2,3}$
represents a universal quadratic lattice $\qf{1,2,3,6}_\Z$ or $\qf{1,2,3,10}_\Z$. Thus the
universality of the Hermitian lattice $\qf{1,2,3}$ over $\Q{-m}$ is proved.

\Case{III-2 a(2)} $m \equiv 3 \pmod{4}$: Then the associated quadratic form of $\qf{1,2,3}$ is
\[
    x_1^2+x_1x_2+\frac{1+m}4x_2^2 + 2\left(y_1^2+y_1y_2+\frac{1+m}4y_2^2\right) +
    3\left(z_1^2+z_1z_2+\frac{1+m}4z_2^2\right).
\]
If we set $x_2=2$ and $y_2=z_2=0$, then the quadratic form becomes $\qf{1,2,3}_\Z+m$. Note that the
quadratic lattice $\qf{1,2,3}_\Z$ represents all positive integers except the form $2^{d} u_5$. If
$n =2^{d} u_5 \ge m$, then $n-m$ is not of the form $2^{d} u_5$ and hence $n-m \to \qf{1,2,3}_\Z$.
It is easily verified that $n \to \qf{1,2,3}$ for $n < m$ by direct calculation. Hence the
Hermitian lattice $\qf{1,2,3}$ is universal. From now through the Case III-2, we may assume that $m
\ne 6, 10; 15, 23, 31, 39$.

\Case{III-2 b} $10 \not\to \qf{1,2,3}$: The next escalation lattice $\ell$ of $\qf{1, 2, 3}$ is of
the form
\[
    \ell \cong
    \begin{pmatrix}
        1 & 0 & 0 & \alpha \\
        0 & 2 & 0 & \beta \\
        0 & 0 & 3 & \gamma \\
        \conj{\alpha} & \conj{\beta} & \conj{\gamma} & 10
    \end{pmatrix}
\]
for some $\alpha, \beta, \gamma$.

\Case{III-2 b(1)} $m \nequiv 3 \pmod 4$: Then the lattice $\ell$
represents a quadratic lattice
\[
    \ell' = \begin{pmatrix}
        1 & 0 & 0 & \Re{\alpha} \\
        0 & 2 & 0 & \Re{\beta} \\
        0 & 0 & 3 & \Re{\gamma} \\
        \Re{\alpha} & \Re{\beta} & \Re{\gamma} & 10
    \end{pmatrix}_\Z.
\]
There are 28 lattices of this type up to isometry and we know that they are all universal by the
fifteen theorem.

\Case{III-2 b(2)} $m \equiv 3 \pmod{4}$: From the condition of positive semi-definiteness, the
escalation lattices $\ell$ of $\qf{1, 2, 3}$ are all inherited for $m \ge 123$. Then since $\ell$
represents the above quadratic lattice  $\ell'$, $\ell$ is universal. We may assume that $m = 35,
43, 47, \dotsc, 119$. Note that $\qf{1,2,3}_\Z$ represents all positive integers except only the
form $2^d u_5$. If $n=2^d u_5 \ge m$, then $n-m$ is not of the form $2^d u_5$. Thus $n \to
\qf{1,2,3}_\Z$. If $n < m$, then we can show that $n$ is represented by each $\ell$ by direct
calculation. Hence the universality of $\ell$ over $\Q{-m}$ is proved.

\Case{III-3} $\qf{1,2,4} \to L$: Note that $\qf{1,2,4}$ represents a quadratic lattice
$\qf{1,2,4}_\Z$ whose truant is $14$.

\Case{III-3 a} $14 \to \qf{1,2,4}$: Then we have $m=6$, $10$, $13$, $14;$ $15$, $23$, $31$, $39$,
$47$, $55$.

\Case{III-3 a(1)} $m \nequiv 3 \pmod{4}$: Since $m = 6,10,13,14$, $\qf{1,2,4}$ over $\Q{-m}$
represents a universal quadratic form $\qf{1,2,4,m}_\Z$. Hence the universality of $\qf{1,2,4}$
over $\Q{-m}$ is proved.

\Case{III-3 a(2)} $m \equiv 3 \pmod{4}$: The associated quadratic form of $\qf{1,2,4}$ is
\[
    x_1^2+x_1x_2+\frac{1+m}4x_2^2 + 2\left(y_1^2+y_1y_2+\frac{1+m}4y_2^2\right) +
    4\left(z_1^2+z_1z_2+\frac{1+m}4z_2^2\right).
\]
If we set $x_2=2$ and $y_2=z_2=0$, then the quadratic form becomes $\qf{1,2,4}_\Z + m$. Note that
the quadratic lattice $\qf{1,2,4}_\Z$ represents all positive integers except the form $2^{d}u_7$.
If $n = 2^{d} u_7 \ge m$, then $n-m$ is represented by $\qf{1,2,4}_\Z$ and hence $n \to
\qf{1,2,4}$. It is easily verified that $\qf{1,2,4}$ represents $n$ for $n < m$. Hence $\qf{1,2,4}$
is universal over $\Q{-m}$. From now through the Case III-3 b, we may assume that $m \ne 6$, $10$,
$13$, $14;$ $15$, $23$, $31$, $39$, $47$, $55$.

\Case{III-3 b} $14 \not\to \qf{1,2,4}$: The next escalation lattice $\ell$ of $\qf{1,2,4}$ is of
the form
\[
    \ell \cong
    \begin{pmatrix}
        1 & 0 & 0 & \alpha \\
        0 & 2 & 0 & \beta \\
        0 & 0 & 4 & \gamma \\
        \conj{\alpha} & \conj{\beta} & \conj{\gamma} & 14
    \end{pmatrix}
\]
for some $\alpha, \beta, \gamma$.

\Case{III-3 b(1)} $m \nequiv 3 \pmod 4$: Then the lattice $\ell$ represents a quadratic lattice
\[
    \ell' = \begin{pmatrix}
        1 & 0 & 0 & \Re{\alpha} \\
        0 & 2 & 0 & \Re{\beta} \\
        0 & 0 & 4 & \Re{\gamma} \\
        \Re{\alpha} & \Re{\beta} & \Re{\gamma} & 14
    \end{pmatrix}_\Z.
\]
There are 54 lattices of this type up to isometry and they are all universal by the fifteen
theorem.

\Case{III-3 b(2)} $m \equiv 3 \pmod{4}$: From the positive semi-definiteness, the escalation
lattices are all inherited for $m \ge 227$. Since $\ell$ represents the above quadratic lattice
$\ell'$, $\ell$ is universal. Thus we may assume that $m = 35, 43, 51, 59, 67,$ $71, 79, 83, 87,
\dotsc, 223$. Note that $\qf{1,2,4}_\Z$ represents all positive integers except only the form
$2^{d} u_7$. We have that $n-m \to \qf{1,2,4}_\Z$ implies $n \to \qf{1,2,4}$. If $n=2^{d} u_7 \ge
m$, then $n-m$ is not of the form $2^{d} u_7$. Thus $n \to \qf{1,2,4}_\Z$. By direct calculation it
can be verified that $n \to \ell$ for $n < m$. Hence the universality of $\ell$ over $\Q{-m}$ is
proved.

\Case{III-4} $\qf{1,2,5} \to L$:  Note that $\qf{1,2,5}$ represents a quadratic lattice
$\qf{1,2,5}_\Z$ whose truant is $10$.

\Case{III-4 a} $10 \to \qf{1,2,5}$: Then we should have $m=6$, $10$; $15$, $23$, $31$, $39$.

\Case{III-4 a(1)} $m \nequiv 3 \pmod{4}$: Since $m = 6, 10$, $\qf{1,2,5}$ over $\Q{-m}$ represents
a quaternary universal quadratic lattice $\qf{1,2,5,m}_\Z$. Hence $\ell$ is universal.

\Case{III-4 a(2)} $m \equiv 3 \pmod{4}$: The associated quadratic form of $\qf{1,2,5}$ is
\[
    x_1^2+x_1x_2+\frac{1+m}4x_2^2 + 2\left(y_1^2+y_1y_2+\frac{1+m}4y_2^2\right) +
    5\left(z_1^2+z_1z_2+\frac{1+m}4z_2^2\right).
\]
If we set $x_2=2$ and $y_2=z_2=0$, then the quadratic form becomes $\qf{1,2,5}_\Z+m$. Similarly
$\qf{1,2,5}_\Z+2m$ and $\qf{1,2,5}_\Z+3m$ are also obtained. The quadratic lattice $\qf{1,2,5}_\Z$
represents all positive integers except the form $5^d u_{-}$. If $m=15$ and $n \ge 3m$, then $n-m$,
$n-2m$, or $n-3m$ is represented by $\qf{1,2,5}_\Z$. If $m=23, 31, 39$ or $n = 5^d u_{-} \geq 2m$,
then $n-m$ or $n-2m$ is represented by $\qf{1,2,5}_\Z$. It is easily verified that $\qf{1,2,5}$
represents $n$ for $n < 3m$ or $n<2m$. Hence $\qf{1,2,5}$ is universal. From now through the Case
III-4 b, we may assume that $m \ne 6, 10; 15, 23, 31, 39$.

\Case{III-4 b} $10 \not\to \qf{1,2,5}$:  The next escalation lattice $\ell$ of $\qf{1,2,5}$ is of
the form
\[  \ell \cong
    \begin{pmatrix}
        1 & 0 & 0 & \alpha \\
        0 & 2 & 0 & \beta \\
        0 & 0 & 5 & \gamma \\
        \conj{\alpha} & \conj{\beta} & \conj{\gamma} & 10
    \end{pmatrix}
\]
for some $\alpha$, $\beta$, $\gamma$.

\Case{III-4 b(1)} $m \nequiv 3 \pmod 4$:  From the positive semi-definiteness, $\ell$ is inherited
if $m \geq 53$. The lattice $\ell$ represents a quadratic lattice
\[
    \begin{pmatrix}
        1 & 0 & 0 & \Re{\alpha} \\
        0 & 2 & 0 & \Re{\beta} \\
        0 & 0 & 5 & \Re{\gamma} \\
        \Re{\alpha} & \Re{\beta} & \Re{\gamma} & 10
    \end{pmatrix}_\Z.
\]
There are 32 quadratic lattices of this type up to isometry. Among these quadratic lattices, 28
quadratic lattices are universal by the fifteen theorem, but the following 4 quadratic lattices are
not. The truants of the following 4 quaternary quadratic lattices are all 15.
\[
\qf{1}_\Z \perp \negthickspace\terlattice{2}{0}{0}{0}{5}{0}{0}{0}{5}_\Z\negthickspace, %
\qf{1}_\Z \perp \negthickspace\terlattice{2}{0}{1}{0}{5}{1}{1}{1}{5}_\Z\negthickspace, %
\qf{1}_\Z \perp \negthickspace\terlattice{2}{0}{1}{0}{5}{2}{1}{2}{8}_\Z\negthickspace, %
\qf{1}_\Z \perp
\negthickspace\terlattice{2}{0}{1}{0}{5}{1}{1}{1}{9}_\Z.
\]
Now assume that $\ell$ represents one of the above quadratic lattices. If $m=13$ or $14$, the
Hermitian lattice $\ell$ represent $15$ and hence they are universal. If $m \ne 13, 14$ and $\ell$
is uninherited, then $\ell$ is one of the followings:

\begin{tabular}{ll}                                       &   \\
$\qf{1} \perp \terlattice20005{\omega_m}0{\comega_m}5$   & if $m = 17, 21, 22,$\\

$\qf{1} \perp \terlattice20105{1\pm\omega_m}1{1\pm\comega_m}5$ & if $m = 17, 21,$\\

$\qf{1} \perp \terlattice20105{2\pm\omega_m}1{2\pm\comega_m}8$ & if $m = 17, 21, 22, 26, 29, 30, 33,$\\

$\qf{1} \perp \terlattice20105{1\pm\omega_m}1{1\pm\comega_m}9$ & if $m = 17, 21, 22, 26, 29, 30, 33, 34, 37, 38, 41.$%

\end{tabular}\\

\noindent Among the above lattices, the only following lattices represent 15. Hence they are
universal by the fifteen theorem.
\begin{align*}
& \qf{1} \perp \terlattice20005{\omega_{22}}0{\comega_{22}}5,
& \qf{1} \perp \terlattice20105{1\pm\omega_{21}}1{1\pm\comega_{21}}5, \\
& \qf{1} \perp \terlattice20105{2\pm\omega_{33}}1{2\pm\comega_{33}}8, %
& \qf{1} \perp
\terlattice20105{1\pm\omega_{41}}1{1\pm\comega_{41}}9.
\end{align*}

\noindent If $\ell$ does not represent 15, then we can obtain a universal \emph{pro forma} quinary
Hermitian lattice by attaching a vector of norm 15 to $\ell$. If $\ell$ is inherited and it is not
universal, then $\ell$ is one of the following 4 lattices whose truants are all 15:
\[
\qf{1} \perp \negthickspace\terlattice{2}{0}{0}{0}{5}{0}{0}{0}{5}, %
\qf{1} \perp \negthickspace\terlattice{2}{0}{1}{0}{5}{1}{1}{1}{5}, %
\qf{1} \perp \negthickspace\terlattice{2}{0}{1}{0}{5}{2}{1}{2}{8}, %
\qf{1} \perp \negthickspace\terlattice{2}{0}{1}{0}{5}{1}{1}{1}{9}.
\]
\noindent In this case, we can obtain a universal lattice by attaching a vector of norm 15. Note
that these quinary universal lattices are inherited if $m \geq 129$.

\Case{III-4 b(2)} $m \equiv 3 \pmod{4}$: From the positive semi-definiteness, the escalation
lattices $\ell$'s are all inherited for $m \ge 203$. We may assume that $m = 35, 43, 47, 51,
\dotsc, 199$. Note that $\qf{1,2,5}_\Z$ represents all positive integers except the form $5^{d}
u_{-}$. If $n=5^{d}u_{-} \ge 3m$, then at least one of $n-m$, $n-2m$ and $n-3m$ is not of the form
$5^{d}u_{-}$. Thus $n \to \qf{1,2,5}$. If $\ell$ represents all positive integers $n < 3m$, then
$\ell$ is universal. If $\ell$ is not universal and uninherited, then $\ell$ is one of the
followings and their conjugates:

\begin{tabular}{ll}
$\qf{1} \perp \terlattice20005{1+\omega}0{1+\comega}8$
                  & if $m = 47,55,151,$ $67 \leq  m \leq 131$,\\
$\qf{1} \perp \terlattice20005{-2+\omega}0{-2+\comega}8$
                  & if $m = 47,55,151,$ $67 \leq  m \leq 131$,\\
$\qf{1} \perp \terlattice20005{2+\omega}0{2+\comega}8$
                  & if $m = 47,55,$ $67 \leq  m \leq 119$,\\
$\qf{1} \perp \terlattice20105{2+\omega}1{2+\comega}9$
                  & if $m = 47,55,$ $67 \leq  m \leq 131$,\\
$\qf{1} \perp \terlattice20005{2+\omega}0{2+\comega}{10}$
                  & if $m = 47,55,$ $67 \leq  m \leq 159.$\\
\end{tabular}

\noindent We can check that $\ell$ represents all positive integers except only 15. Hence we can
obtain a universal lattice by attaching a vector of norm $15$ to $\ell$. If $\ell$ is inherited and
$\ell$ is not universal, then we also obtain a universal lattice by the same process.

\Case{III-5} $\qf{1}\perp\binlattice2114 \to L$: Note that $\qf{1}\perp\binlattice2114$ represents
a quadratic lattice $\qf{1}_\Z\perp\binlattice2114_\Z$ whose truant is $7$.

\Case{III-5 a} $7 \to \qf{1}\perp\binlattice2114$: Then we only have $m=6$. Since
$\qf{1}\perp\binlattice2114$ represents a universal quadratic lattice $\qf{1,6}_\Z \perp
\binlattice2114_\Z$, $\qf{1}\perp\binlattice2114$ is universal. From now through the Case III-5, we
may assume that $m \ne 6$.

\Case{III-5 b} $7 \not\to \qf{1}\perp\binlattice2114$: The next escalation lattice $\ell$ of
$\qf{1}\perp\binlattice2114$ is of the form
\[ \ell \cong
    \begin{pmatrix}
        1 & 0 & 0 & \alpha \\
        0 & 2 & 1 & \beta \\
        0 & 1 & 4 & \gamma \\
        \conj{\alpha} & \conj{\beta} & \conj{\gamma} & 7
    \end{pmatrix}
\]
for some $\alpha, \beta$ and $\gamma$.

\Case{III-5 b(1)} $m \nequiv 3 \pmod{4}$: From the positive semi-definiteness, $\ell$ is inherited
if $m \ge 29$. And $\ell$ represents a quadratic lattice
\[
    \begin{pmatrix}
        1 & 0 & 0 & \Re{\alpha} \\
        0 & 2 & 1 & \Re{\beta} \\
        0 & 1 & 4 & \Re{\gamma} \\
        \Re{\alpha} & \Re{\beta} & \Re{\gamma} & 7
    \end{pmatrix}_\Z.
\]
There are 30 quadratic lattices of this type up to isometry. These quadratic lattices are universal
except the following lattice
\[
    \ell' = \qf{1}_\Z \perp \terlattice{2}{1}{0}{1}{4}{1}{0}{1}{5}_\Z.
\]
Note that the truant of $\ell'$ is 10 and $\ell'$ represents all numbers 1 to 15 except 10. Assume
that $\ell$ represents $\ell'$. If $m=10$, then the Hermitian lattice $\ell$ represents 10, hence
$\ell$ is a quaternary universal Hermitian lattice. From the positive semi-definiteness, if $m \geq
17$, then $\ell$ is inherited. If $m = 13, 14$ and $\ell$ is uninherited, then
\[
    \ell \cong \qf{1} \perp \terlattice21014{1\pm\omega_m}0{1\pm\comega_m}5.
\]
Since $\ell$ represents all positive integers 1 through 15 except only the truant 10 of $\ell$, the
next escalation lattice of $\ell$ is a \emph{pro forma} quinary universal Hermitian lattice which
can be obtained by attaching a vector of norm 10. In this case, the universal Hermitian lattice is
inherited if $m \ge 53$. If $\ell$ is inherited and it is not universal, then $\ell$ is the
following lattice whose truant is 10:
\[
    \qf{1} \perp \terlattice{2}{1}{0}{1}{4}{1}{0}{1}{5}.
\]
Hence we can obtain a universal lattice via attaching a vector of norm 10.

\Case{III-5 b(2)} $m \equiv 3 \pmod{4}$: From the positive semi-definiteness, the escalation
lattices $\ell$'s are all inherited for $m \ge 115$. Thus we may assume that $m = 23, 31, 35,
\dotsc, 111$. Note that $\qf{1}_\Z \perp \binlattice2114_\Z$ represents positive integers $n$ when
$n \ne 7^{d} u_{-}$ and $n \nequiv 7,10 \pmod{12}$. The  lattice $\qf{1} \perp \binlattice2114$
represents $\qf{1}_\Z \perp \binlattice2114_\Z + km$ with $k=1$, $2$, $3$. If $n$ is not
represented by $\qf{1}_\Z \perp \binlattice2114_\Z$, then at least one of $n-m$, $n-2m$ and $n-3m$
is represented by $\qf{1}_\Z \perp \binlattice2114_\Z$ for $n \ge 3m$. It can be verified that
$\ell$ represents all positive integers smaller than $3m$. Hence the universality of $\ell$ is
proved. If $m \geq 115$, then all $\ell$'s are universal except $\qf{1} \perp
\terlattice{2}{1}{0}{1}{4}{1}{0}{1}{5}$. We can obtain a universal \emph{pro forma} quinary lattice
by attaching a vector of norm 10 to this exception.

\Case{III-6} $\qf{1}\perp\binlattice2115 \to L$: Note that $\qf{1}\perp\binlattice2115$ represents
a quadratic lattice $\qf{1}_\Z\perp\binlattice2115_\Z$ whose truant is $7$.

\Case{III-6 a} $7 \to \qf{1}\perp\binlattice2115$: We have only $m=6$. Since
$\qf{1}\perp\binlattice2115$ represents a universal quadratic lattice $\qf{1,6}_\Z \perp
\binlattice2115_\Z$, $\qf{1}\perp\binlattice2115$ is universal. From now through the Case III-6, we
may assume that $m \ne 6$.

\Case{III-6 b} $7 \not\to \qf{1}\perp\binlattice2115$: The next
escalation lattice $\ell$ of $\qf{1}\perp\binlattice2115$ is of
the form
\[
    \ell \cong
    \begin{pmatrix}
        1 & 0 & 0 & \alpha \\
        0 & 2 & 1 & \beta \\
        0 & 1 & 5 & \gamma \\
        \conj{\alpha} & \conj{\beta} & \conj{\gamma} & 7
    \end{pmatrix}
\]
for some $\alpha, \beta$, $\gamma$.

\Case{III-6 b(1)}  $m \nequiv 3 \pmod{4}$: From the positive
semi-definiteness, $\ell$ is inherited if $m \ge 37$. And $\ell$
represents a quadratic lattice
\[
    \begin{pmatrix}
        1 & 0 & 0 & \Re{\alpha} \\
        0 & 2 & 1 & \Re{\beta} \\
        0 & 1 & 5 & \Re{\gamma} \\
        \Re{\alpha} & \Re{\beta} & \Re{\gamma} & 7
    \end{pmatrix}_\Z.
\]
There are 16 quadratic lattices of this type up to isometry. These
quadratic lattices are universal except the following lattice
\[
    \ell' = \qf{1}_\Z \perp \terlattice{2}{1}{0}{1}{5}{1}{0}{1}{5}_\Z.
\]
Note that the truant of $\ell'$ is 15 and $\ell'$ represents all numbers 1 to 14. Assume $\ell$
represents $\ell'$. If $\ell$ is uninherited, then $m = 10, 13, 14, 17, 21$ and
\[
    \ell \cong \qf{1} \perp \terlattice21015{1\pm\omega_m}0{1\pm\comega_m}5.
\]
Note that if $m= 10, 13, 14, 21$, then $\ell$ represents 15. Hence $\ell$ is a quaternary universal
lattice. If $m = 17$, then it is isometric to the lattice in the Case III-4 b(1). If $\ell$ is
inherited, and it is not universal, then $\ell$ is the following lattice whose truant is 15:
\[
    \qf{1} \perp \terlattice{2}{1}{0}{1}{5}{1}{0}{1}{5}.
\]
Thus we can obtain a universal lattice via attaching a vector of norm 15.

\Case{III-6 b(2)} $m \equiv 3 \pmod{4}$: From the positive semi-definiteness, the escalation
lattices $\ell$'s are inherited for $m \ge 147$. Thus we may assume that $m = 23, 31, 35, \dotsc,
139$. Note that $\qf{1}_\Z \perp \binlattice2115_\Z$ represents all positive integers $n$ except
the form $2^{e} u_7$. The associated (nonclassical) quadratic lattice of $\qf{1} \perp
\binlattice2115$ represents $\qf{1}_\Z \perp \binlattice2115_\Z \perp \qf{m}_\Z$. If $n=2^{e} u_7 >
4m$ is not represented by $\qf{1}_\Z \perp \binlattice2115_\Z$, then at least one of $n-m$ and $n -
4m$ is represented by $\qf{1}_\Z \perp \binlattice2115_\Z$. It is also verified that each $\ell$
represents all positive integers smaller than $4m$. Hence the universality of $\ell$'s is proved.
If $m \ge 147$, then all $\ell$'s are universal except $\qf{1} \perp
\terlattice{2}{1}{0}{1}{5}{1}{0}{1}{5}$. We can obtain a universal \emph{pro forma} quinary lattice
by attaching a vector of norm 15 to this exception.

\Case{III-7} Now we investigate the lattices
\[
    \qf{1} \perp \binlattice2{\omega}{\comega}b   \text{ and } \qf{1} \perp \binlattice2{-1+\omega}{-1+\comega}b
\]
over $\Q{-m}$ in \ref{tbl:escalation_lattices_in_Case_III}.  It is known that all the binary
Hermitian lattices in Table \ref{tbl:escalation_lattices_in_Case_III} are universal (see
\cite{hI_00}) except for the last two over $\Q{-39}$ and they are listed in the end of this paper
(see Table \ref{tbl:binary_universal_Hermitian_lattices}). If $m \nequiv 3 \pmod{4}$, i.e.,
$m=6,10$, then the universalities of ternary lattices in Table
\ref{tbl:escalation_lattices_in_Case_III} are checked by the fifteen theorem. Now assume that $m
\equiv 3 \pmod{4}$, i.e., $m=15,23,31,35,39$. Since
\[
    \qf{1} \perp \binlattice2{-1+\omega}{-1+\comega}b
    = \conj{\qf{1} \perp \binlattice2{\omega}{\comega}b},
\]
it is enough to check the universality of $\ell \cong \qf{1} \perp \binlattice2{\omega}{\comega}b$.
Each associated quadratic lattice of $\ell$ has a sublattice $\ell'$ of class number $1$ (see Table
\ref{tbl:Sublattices of class number 1_in_Case_III-7}). The universality of each lattice $\ell'$ is
proved by the method similar to the previous case.

\begin{table}[!h]
\small
\begin{tabular}{l|l|l}
    \hline
    \multicolumn{1}{c|}{\normalsize field}
        & \multicolumn{1}{c|}{\normalsize lattice}
        & \multicolumn{1}{c}{\normalsize subtrahend} \\ \hline
    {\normalsize $\Q{-15}$}
        & $\terlattice10002{\omega}0{\comega}3 \leftarrow \terlattice10002404{12}_\Z \negthickspace\cong
        \terlattice100020004_\Z \negthickspace\leftarrow n \ne 2^du_7$
        & $15$ \\
    ~
        & $\terlattice10002{\omega}0{\comega}4 \leftarrow \terlattice10002404{16}_\Z \negthickspace\cong
        \terlattice100020008_\Z \negthickspace\leftarrow n \ne u_5, 2^eu_7$
        & $15, 2\cdot15$ \\
    ~
        & $\terlattice10002{\omega}0{\comega}5 \leftarrow \terlattice100084045_\Z \negthickspace\cong
        \terlattice100051015_\Z \negthickspace\leftarrow n \ne 2^du_5$
        & $15$ \\
    \hline
    {\normalsize $\Q{-23}$}
        & $\terlattice10002{\omega}0{\comega}4 \leftarrow \terlattice10002606{24}_\Z \negthickspace\cong
        \terlattice100020006_\Z \negthickspace\leftarrow n \ne 2^eu_5$
        & $23$ \\
    ~
        & $\terlattice10002{\omega}0{\comega}5 \leftarrow \terlattice1000{12}6065_\Z \negthickspace\cong
        \terlattice100051015_\Z \negthickspace\leftarrow n \ne 2^du_5$
        & $23$ \\
    \hline
    {\normalsize $\Q{-31}$}
        & $\terlattice10002{\omega}0{\comega}5 \leftarrow \terlattice10002808{40}_\Z \negthickspace\cong
        \terlattice100020008_\Z \negthickspace\leftarrow n \ne u_5, 2^eu_7$
        & $15, 2\cdot15$ \\
    \hline
    {\normalsize $\Q{-35}$}
        & $\terlattice10002{\omega}0{\comega}5 \leftarrow \terlattice10002909{45}_\Z \negthickspace\cong
        \terlattice100021015_\Z \negthickspace\leftarrow n \ne 2^eu_7$
        & $35, 2^2\cdot35$ \\
    \hline
\end{tabular}
\vspace{1ex}
\caption{Sublattices of class number 1}%
\label{tbl:Sublattices of class number 1_in_Case_III-7}%
\end{table}

The last exceptional lattice is
\[
    \ell \cong \qf{1} \perp\binlattice2{\omega}{\comega}5
\]
over $\Q{-39}$. The associated quadratic form of this lattice is
\[
    N: x_1^2 + x_1 x_2 + 10 x_2^2 + 2y_1^2 + y_1 z_1 + 5z_1^2.
\]
Since $2N \to N$, we only need to show that $N$ represents all odd positive integers. If we set
$x_2 = 2t$ for some $t \in \Z$, then $N$ represents
\[
    \qf{1}_\Z \perp \binlattice2{1/2}{1/2}5_\Z \perp \qf{39}_\Z.
\]
For sufficiently large $n$, we will show that
\[
    n - 39 s^2 \to  N' = \qf{1}_\Z \perp
    \binlattice2{1/2}{1/2}5_\Z
\]
for suitable $s$. It implies the desired result: $n \to \ell$ over $\Q{-39}$. Note that $N'$
represents two quadratic sublattices which are in the same genus by Brandt-Intrau table
\cite{hB_oI}:
\[
    \terlattice10008202{20}_\Z, ~~\terlattice402051219_\Z .
\]
Hence if we show that $n - 39s^2$ is represented by the genus, then we can say $n \to \ell$. We
have that the genus represents all positive integers $n$ when $n \equiv 0,1 \pmod{4}$ and $n \ne
13^d u_{+}$. For other odd integers $n \ge 39 \cdot 6^2$, here are choices for $s$ as follows:

\begin{table}[!h]
\begin{tabular}{l|l|l} \hline%
\multicolumn{1}{c|}{$n$} & \multicolumn{1}{c|}{$u_{\pm}$} & subtrahend \\ \hline%
$n \equiv 1 \pmod{4}, n=13^d u_{+}$ &
    $u_{+} \equiv 1,4,10 \pmod{13}$ & $39\cdot 2^2$ \\
~ &
    $u_{+} \equiv 3 \pmod{13}$ & $39\cdot 4^2$ \\
~ &
    $u_{+} \equiv 9,12 \pmod{13}$ & $39\cdot 6^2$ \\ \hline
$n \equiv 3 \pmod{4}, n=13^e u_{+}$ &
    $u_{+} \equiv 1,4,9,10,12 \pmod{13}$ & $39\cdot 1^2$ \\
~ &
    $u_{+} \equiv 3 \pmod{13}$ & $39\cdot 3^2$ \\
    \hline
$n \equiv 3 \pmod{4}, n=13^e u_{-}$ &
    $u_{-} \equiv 2,5,6,7,8,11 \pmod{13}$ & $39\cdot 1^2$ \\
    \hline
$n \equiv 3 \pmod{4}, n=13^d u_{+}$ &
    $u_{+} \equiv 1,9,10 \pmod{13}$ & $39\cdot 1^2$ \\
~ &
    $u_{+} \equiv 3,12 \pmod{13}$ & $39\cdot 3^2$ \\
~ &
    $u_{+} \equiv 4 \pmod{13}$ & $39\cdot 5^2$ \\ \hline
$n \equiv 3 \pmod{4}, n=13^d u_{-}$ &
    $u_{-} \equiv 2 \pmod{13}$ & $39\cdot 5^2$ \\
~ &
    $u_{-} \equiv 5,8,11 \pmod{13}$ & $39\cdot 1^2$ \\
~ &
    $u_{-} \equiv 6,7 \pmod{13}$ & $39\cdot 3^2$ \\ \hline
\end{tabular}
\end{table}

For $n < 39 \cdot 6^2$, we can check that $\ell$ represents all integers $n$ except $13$ and $91$
by direct calculation. If we attach a vector of norm $13$ to $\ell$, the escalation lattices also
represent $91$. Hence they are universal Hermitian lattices.

\begin{Rmk}
In this article, we computed all new universal Hermitian lattices which does not contain other
universal lattices. For example, over $\Q{-5}$, there are $2$ binary universal Hermitian lattices
as followings, up to isometry \cite{hI_00}:
\[
    \qf{1,2},~~
    \qf{1} \perp \binlattice2{-1+\omega}{-1+\conj{\omega}}3.
\]
The first lattice is found in the procedure of escalation in this article, but the second one is
excluded because it contain a universal Hermitian lattice $\qf{1,2}$ as a sublattice.
\end{Rmk}

\begin{Rmk}
We computed all quaternary escalation lattices and corresponding truants (See Table \ref{tbl:Proof
of universality of escalator}). We can obtain \emph{pro forma} quinary universal Hermitian lattices
by attaching a vector whose norm is its truant, to each quaternary escalation lattice. But these
are too numerous to list all.
\end{Rmk}

From the above results we obtain the following theorems.

\begin{Thm}
For all positive square-free integers $m$, if a Hermitian lattice over $\Q{-m}$ represents the
following critical numbers, then it is universal.\\

\begin{center}
\begin{tabular}{l|l}  \hline
critical numbers                & $m$ \\ \hline %
$1, 2$                          & $3, 11$ \\
$1, 3$                          & $1, 7$  \\
$1, 5$                          & $2$ \\
$1, 2, 3$                       & $5, 19$ \\
$1, 2, 3, 5$                    & $6$ \\
$1, 2, 3, 5, 7$                 & $15, 23$ \\
$1, 2, 3, 5, 6, 7$              & $10, 31$ \\
$1, 2, 3, 5, 6, 7, 10$          & $13, 14$ \\
$1, 2, 3, 5, 6, 7, 13$          & $39$\\
$1, 2, 3, 5, 6, 7, 10, 14$      & $35, 43, 51, 59$\\
$1, 2, 3, 5, 6, 7, 10, 15$      & $55$\\
$1, 2, 3, 5, 6, 7, 10, 14, 15$  & otherwise \\
\hline
\end{tabular}
\end{center}
\end{Thm}

From this theorem, we have the criterion on the universality of Hermitian lattices.

\begin{Cor}[The fifteen theorem for universal Hermitian lattices]
If a positive definite Hermitian lattice represents $1$, $2$, $3$, $5$, $6$, $7$, $10$, $13$, $14$,
and $15$, then it is universal.
\end{Cor}

Unlike the fifteen theorem for universal \emph{quadratic} forms, this theorem involves the number
$13$ for the case of $\Q{-39}$.

\begin{Thm}
If we set
\[
    u_{m} := \min\left\{ \rank{L}: L \text{ is a universal Hermitian lattice over }\Q{-m} \right\}
\]
for each positive square-free integer $m$, then

\[
\begin{array}{c|llllllllllll}  \hline %
u_{m} &       &     &      &      &      & m    &      &      &      &      &      &       \\
\hline
  2   & 1,    & 2,  & 3,   & 5,   & 6,   & 7,   & 10,  & 11,  & 15,  & 19,  & 23,  & 31,   \\
\hline
  3   & 13,   & 14, & 17,  & 21,  & 22,  & 26,  & 29,  & 30,  & 34,  & 35,  & 39,  & 41,  \\
      & 43,   & 46, & 47,  & 51,  & 55,  & 59,  & 71,  & 79,  & 83,  & 87,  & 91,  & 95,  \\
      & 103,  & 107,& 111, & 115, & 119, & 127, & 131, & 135, & 139, & 143, & 147, & 151, \\
      & 155,  & 159,& 167, & 171, & 175, & 179, & 183, & 187, & 191, & 199, & 207, & 215, \\
      & 223,  \\
\hline
  4   & 33,   & 37,  & 38,  & 42,  & 53,  & 57,  & 58,   & 61,  & 62,  & 65,  & 66,  & 67,  \\
      & 69,   & 70,  & 73,  & 74,  & 77,  & 78,  & 82,   & 85,  & 86,  & 89,  & 93,  & 94,  \\
      & 97,   & 101, & 102, & 105, & 106, & 109, & 110,  & 113, & 114, & 118, & 122, & 123, \\
      & 129, & 130, & 133, & 134, & 137, & 138,  & 141,  & 142, & 145, & 146, & 149, & 154, \\
      & 157, & 158, & 161, & 163, & 165, & 166,  & 170,  & 173, & 174, & 177, & 178, & 181, \\
      & 182, & 185, & 186, & 190, & 193, & 194,  & 195,  & 197, & 201, & 202, & 203, & 205, \\
      & 206, & 209, & 210, & 211, & 213, & 214,  & 217,  & 218, & 219, & 221, & 222, & 226, \\
      & 227, & 229, & 230, & 231  & \text { or } &m \ge &235.      \\ \hline
\end{array}
\]

\end{Thm}
Note that $u_{m} \leq 4$, since $\qf{1,1,1,1}$ is an inherited universal Hermitian lattice over
$\Q{-m}$ for all positive square-free integers $m$.


\begin{table}[p]
\begin{tabular}{lll} \hline
$\Q{-m}$          & binary universal lattices \\ \hline%
\mystrut$\Q{-1}$  & $\qf{1,1}$, $\qf{1,2}$, $\qf{1,3}$ \\
\mystrut$\Q{-2}$  & $\qf{1,1}$, $\qf{1,2}$, $\qf{1,3}$, $\qf{1,4}$, $\qf{1,5}$ \\
\mystrut$\Q{-3}$  & $\qf{1,1}$, $\qf{1,2}$ \\
\mystrut$\Q{-5}$  & $\qf{1,2}$, $\qf{1}\perp\binlattice{2}{-1+\omega_{5}}{-1+\comega_{5}}3$ \\
\mystrut$\Q{-6}$  & $\qf{1}\perp\binlattice{2}{\omega_6}{\comega_6}3$ \\
\mystrut$\Q{-7}$  & $\qf{1,1}$, $\qf{1,2}$, $\qf{1,3}$ \\
\mystrut$\Q{-10}$ & $\qf{1}\perp\binlattice{2}{\omega_{10}}{\comega_{10}}5$ \\
\mystrut$\Q{-11}$ & $\qf{1,1}$, $\qf{1,2}$ \\
\mystrut$\Q{-15}$ & $\qf{1}\perp\binlattice{2}{\omega_{15}}{\comega_{15}}2$ \\
\mystrut$\Q{-19}$ & $\qf{1,2}$ \\
\mystrut$\Q{-23}$ &
$\qf{1}\perp\binlattice{2}{\omega_{23}}{\comega_{23}}3$,
$\qf{1}\perp\binlattice{2}{-1+\omega_{23}}{-1+\comega_{23}}3$ \\
\mystrut$\Q{-31}$ &
$\qf{1}\perp\binlattice{2}{\omega_{31}}{\comega_{31}}4$,
$\qf{1}\perp\binlattice{2}{-1+\omega_{31}}{-1+\comega_{31}}4$
\\ \hline
\end{tabular} \vspace{2ex}
\caption{Binary universal Hermitian lattice}%
\label{tbl:binary_universal_Hermitian_lattices}
\end{table}

\begin{scriptsize}
\begin{table}[p]
\begin{tabular}{lcl}  \hline
Escalation lattice & Truant & $m$          \\
\hline
$\qf{1}$          & $2$      & if  $m \ne$  $1,2;7,$\\
                  & $3$      & if  $m =$    $1;7,$\\
                  & $5$      & if  $m =$    $2,$\\

$\qf{1,1}$        & $3$      & if  $m \ne$  $1,2;3,7,11,$\\

$\qf{1,1,1}$      & $7$      & if  $m \ne$  1,2,5,6;3,7,11,15,19,23,\\
$\qf{1,1,2}$      & $14$     & if  $m \ne$  1,2,5,6,10,13,14;3,7,11,15,19,23,31,35,39,43,47,51,55,\\

$\qf{1,1,3}$      & $6$      & if  $m \ne$  1,2,5,6;3,7,11,15,19,23,\\

$\qf{1,2}$        & $5$      & if  $m \ne$  1,2,5;3,7,11,19,\\

$\qf{1,2,2}$      & $7$      & if  $m \ne$  1,2,5,6;3,7,11,15,19,\\
$\qf{1,2,3}$      & $10$     & if  $m \ne$  1,2,5,6,10;3,7,11,15,19,23,31,39,\\
$\qf{1,2,4}$      & $14$     & if  $m \ne$  1,2,5,6,10,13,14;3,7,11,15,19,23,31,39,47,55,\\
$\qf{1,2,5}$      & $10$     & if  $m \ne$  1,2,5,6,10;3,7,11,15,19,23,31,39,\\

$\qf{1} \perp \terlattice200050005$
                  & 15     & \parbox{20em}{
                             if  $m \equiv 1,2 \pmod{4}$ and $m \geq 22$,\\
                             if  $m \equiv 3 \pmod{4}$ and $m = 47,55,$ $m \geq 67$,  \\
                             } \\
$\qf{1} \perp \terlattice201051115$
                  & 15     & \parbox{20em}{
                             if $m \equiv 1,2 \pmod{4}$ and $m \geq 21$,\\
                             if $m \equiv 3 \pmod{4}$ and $m = 47,55,$ $m \geq 67$,  \\
                             } \\
$\qf{1} \perp \terlattice201052128$
                  & 15     & \parbox{20em}{
                             if $m \equiv 1,2 \pmod{4}$ and $m \geq 33$,\\
                             if $m \equiv 3 \pmod{4}$ and $m = 47,55,$ $m \geq 67$, \\
                             } \\
$\qf{1} \perp \terlattice201051119$
                  & 15     & \parbox{20em}{
                             if $m \equiv 1,2 \pmod{4}$ and $m \geq 41$,\\
                             if $m \equiv 3 \pmod{4}$ and $m = 47,55,$ $m \geq 67$, \\
                             } \\
$\qf{1} \perp \terlattice20005{\pm\omega}0{\pm\comega}5$
                  & 15     & if  $m =$    17, 21, \\
$\qf{1} \perp \terlattice20105{1\pm\omega}1{1\pm\comega}5$
                  & 15     & if $m =$    17, \\
$\qf{1} \perp \terlattice20105{2\pm\omega}1{2\pm\comega}8$
                  & 15     & if $m =$    17, 21, 22,26,29,30,\\
$\qf{1} \perp \terlattice20105{1\pm\omega}1{1\pm\comega}9$
                  & 15     & if $m =$    17,21,22,26,29,30,33,34,37,38,\\
$\qf{1} \perp \terlattice20005{1+\omega}0{1+\comega}8$
                  & 15     & if $m \equiv 3 \pmod{4}$ and $m = 47,55,151,$ $67 \leq  m \leq 131$,\\
$\qf{1} \perp \terlattice20005{-2+\omega}0{-2+\comega}8$
                  & 15     & if $m \equiv 3 \pmod{4}$ and $m = 47,55,151,$ $67 \leq  m \leq 131$,\\
$\qf{1} \perp \terlattice20005{2+\omega}0{2+\comega}8$
                  & 15     & if $m \equiv 3 \pmod{4}$ and $m = 47,55,$ $67 \leq  m \leq 119$,\\
$\qf{1} \perp \terlattice20105{2+\omega}1{2+\comega}9$
                  & 15     & if $m \equiv 3 \pmod{4}$ and $m = 47,55,$ $67 \leq  m \leq 131$,\\
$\qf{1} \perp \terlattice20005{2+\omega}0{2+\comega}{10}$
                  & 15     & if $m \equiv 3 \pmod{4}$ and $m = 47,55,$ $67 \leq  m \leq 159,$\\
$\qf{1} \perp \binlattice2114$
                  & 7      & if $m \ne$  1,2,5,6,10;3,7,11,19,\\

$\qf{1} \perp \terlattice210141015$
                  & 10     & \parbox{20em}{
                             if $m \equiv 1,2 \pmod{4}$ and $m \geq 17$,\\
                             if $m \equiv 3 \pmod{4}$ and $m \geq 115$, \\
                             } \\
$\qf{1} \perp \terlattice21014{1\pm\omega}0{1\pm\comega}5$
                  & 10     & if $m =$     13,14,\\

$\qf{1} \perp \binlattice2115$
                  & 7      & if $m \ne$  1,2,5,6;3,7,11,19,\\

$\qf{1} \perp \terlattice210151015$
                  & 15     & \parbox{20em}{
                             if $m \equiv 1,2 \pmod{4}$ and $m \geq 21$,\\
                             if $m \equiv 3 \pmod{4}$ and $m \geq 147$, \\
                             } \\
$\qf{1} \perp \terlattice21015{1\pm\omega}0{1\pm\comega}5$
                  & 15     & if $m =$     17, \\

$\qf{1} \perp \binlattice2{\omega}{\comega}5$
                  & 13     & if $m =$ 39.\\ \hline
\end{tabular} \vspace{2ex}
\caption{Truants of escalation lattices}
\label{tbl:Proof of universality of escalator}
\end{table}
\end{scriptsize}

\end{document}